\documentclass{SCAEOL}
\numberwithin{equation}{section}
\usepackage{amsfonts}
\usepackage{enumerate}
\usepackage{cases}
\usepackage{type1cm}

\def\mr{\mathbb{R}}
\def\mrn{\mathbb{R}^n}
\def\mrm{\mathbb{R}^m}
\def\mn{\mathbb{N}}
\def\mrd{\mathbb{R}^d}

\def\mmf{\mathbb{F}}
\def\mf{\mathcal{F}}

\def\mb{\mathcal{B}}

\def\mh{\mathbb{H}}
\def\mS{\mathrm{S}}
\def\mT{\mathrm{T}}
\def\ms{\mathbb{S}}

\def\mss{\widetilde{\mathbb{S}}}

\def\mt{\mathbb{T}}

\def\md{\mathbb{D}}
\def\ml{\mathbb{L}}
\def\me{\mathbb{E}}
\def\mmu{\mathcal{U}}

\def\t{\tau}
\def\mx{\mathcal{X}}
\def\dd{\mathcal{D}}
\def\ephs{\varepsilon}

\DeclareMathOperator*{\esssup}{ess~sup}

\newcommand{\inner}[2]{\left\langle#1,#2\right\rangle}

\begin{document}

\Year{20xx} %
\Month{X}
\Vol{xx} %
\No{xx} %
\BeginPage{1} %
\EndPage{xx} %
\AuthorMark{Zhang H {\it et al.}}
\ReceivedDay{xx, 20xx}
\AcceptedDay{xx, 20xx}

\title{Some results on pointwise second-order necessary conditions for stochastic optimal controls}


\author[1]{ZHANG Haisen}{}
\author[2]{ZHANG Xu}{}

%
\address[{\rm1}]{School of Mathematics, Sichuan University, Chengdu {\rm 610064}, China;}
\address[{\rm2}]{Yangtze Center of Mathematics, Sichuan University, Chengdu {\rm 610064}, China;}
\Emails{haisenzhang@yeah.net, zhang\_xu@scu.edu.cn}

\maketitle


{\begin{center}
\parbox{14.5cm}
{\begin{abstract}
The purpose of this paper is to derive some pointwise second-order necessary conditions for stochastic optimal controls in the general case that the control variable enters into both the drift and the diffusion terms. When the control region is convex, a pointwise second-order necessary  condition for stochastic singular optimal controls in the classical sense is established; while when the control region is allowed to be nonconvex, we obtain a pointwise second-order necessary condition for stochastic singular optimal controls in the sense of Pontryagin-type maximum principle. It is found that, quite different from the first-order necessary conditions, the correction part of the solution to the second-order adjoint equation appears in the pointwise second-order necessary conditions whenever the diffusion term depends on the control variable, even if the control region is convex.
\vspace{-3mm}
\end{abstract}}
\end{center}}

\keywords{stochastic optimal control, needle variation, Pontryagin-type maximum principle, pointwise second-order necessary condition, Malliavin calculus.}

\MSC{Primary 93E20; Secondary 60H07, 60H10.}

\renewcommand{\baselinestretch}{1.2}
\begin{center} \renewcommand{\arraystretch}{1.5}
{\begin{tabular}{lp{0.8\textwidth}} \hline \scriptsize
{\bf Citation:}\!\!\!\!&\scriptsize Zhang H, Zhang X.  Some results on pointwise second-order necessary conditions for stochastic optimal controls. Sci China Math, 20xx, xx, doi: 10.1007/s11425-000-0000-0\vspace{1mm}
\\
\hline
\end{tabular}}\end{center}

\baselineskip 11pt\parindent=10.8pt  \wuhao
\section{Introduction}
Let $T>0$ be a fixed constant and $(\Omega,\mf, \mmf, P)$  be a complete filtered
probability space (satisfying the usual conditions), on which a $1$-dimensional standard Wiener process $W(\cdot)$ is defined such that $\mmf=\{\mf_{t} \}_{0\le t\le T}$ is the natural filtration generated by $W(\cdot)$ (augmented by all $P$-null sets in $\mf$).
We consider the following controlled stochastic differential equation:
\begin{equation}\label{controlsys}
\left\{
\begin{array}{l}
dx(t)=b(t,x(t),u(t))dt+\sigma(t,x(t),u(t))dW(t),\ \ \ t\in[0,T],\\
x(0)=x_0,
\end{array}\right.
\end{equation}
with a cost functional
\begin{equation}\label{costfunction}
J(u(\cdot))=\me\Big[\int_{0}^{T}f(t,x(t),u(t))dt+h(x(T))\Big],
\end{equation}
where the stochastic process $u(\cdot)$ is the control valued in a region $U\subset \mrm$ ($m\in \mn$), the stochastic process $x(\cdot)$ is the state valued in $\mr^n$ ($n\in\mn$), and $b,\ \sigma:[0,T]\times \mrn\times U\to \mrn$, $f:[0,T]\times \mrn\times U\to \mr$ and $h:\mrn\to \mr$ are given functions satisfying suitable conditions (to be specified later).

Any adapted stochastic process valued in $U$ is called an admissible control and we denote by $\mmu_{ad}$ the set of admissible controls. Under some standard assumptions, the corresponding state $x(\cdot)$ (of (\ref{controlsys})) is uniquely defined by any given initial datum $x_{0}\in\mr^n$ and admissible control $u(\cdot)\in \mmu_{ad}$.

The stochastic optimal control problem considered in this paper is to find a control
$\bar{u}(\cdot)\in\mmu_{ad}$ such that
\begin{equation}\label{minimum J}
J(\bar{u}(\cdot))=\inf_{u(\cdot)\in \mmu_{ad}}J(u(\cdot)).
\end{equation}
Any $\bar{u}(\cdot)\in \mmu_{ad}$ satisfying (\ref{minimum J}) is called an optimal control. The corresponding state $\bar{x}(\cdot)$ is called an optimal sate, and $(\bar{x}(\cdot),\bar{u}(\cdot))$ is called an optimal pair (for the above optimal control problem).

In the optimal control theory, one of the central issues is to establish necessary conditions for optimal controls. Some early studies on the first-order necessary condition for stochastic optimal controls in the case that the diffusion term is independent of the control variable can be found in \cite{Haussmann76,Kushner72}. As for early works on the same problem but in the case of the diffusion term containing the control variable, we refer to \cite{Bensoussan81,Bismut78}.

Compared to the deterministic setting, new phenomenon appears when deriving the first-order necessary condition (for stochastic optimal controls) for the case that the diffusion term contains the control variable and the control region is possibly nonconvex. Indeed, for the case of nonconvex control region, the needle variation, which is essentially a perturbation technique with respect to the measure, has to be used as a variation of the optimal controls. When the control variable appears in the diffusion term, the state variation is only an infinitesimal of order $\frac{1}{2}$ with respect to the perturbation measure $\varepsilon$ (as $\varepsilon\to 0^+$).
Therefore, to establish the first-order necessary condition (with respect to this perturbation measure), one needs to expand the cost functional up to order two, and therefore two adjoint equations need to be introduced. The first-order necessary condition for this general case, called general (stochastic) Pontryagin-type maximum principle, was obtained by Peng \cite{Peng90}.

Similar to the deterministic case, the first-order necessary condition provides a basic tool to study the properties of stochastic optimal controls and solve them numerically. However, exactly as its deterministic counterpart or even as that in the classical calculus, in some situation the first-order necessary condition for stochastic optimal controls may be trivial and therefore it cannot provide enough information to find the desired optimal controls. Consequently, it is quite natural to study the second-order necessary condition for stochastic optimal controls, especially after the first-order necessary condition for the general situation was established in \cite{Peng90}. Unfortunately, to the authors' best knowledge, very few works are available in this respect.

Recently, Tang \cite{Tang10} derived a pointwise second-order necessary condition for stochastic optimal controls with nonconvex control regions for the special case that the diffusion term is independent of the control variable; while Bonnans and Silva \cite{Bonnans12} obtained an integral-type (rather than the more desired pointwise-type) second-order necessary condition for stochastic optimal controls with the control variable entering into the diffusion terms, but the control region is assumed to be convex. It seems us that \cite{Bonnans12} and \cite{Tang10} are the only two publications on second-order necessary conditions for stochastic optimal controls. Also, as far as we know, there is no article addressed to the \emph{pointwise} second-order necessary condition for stochastic optimal controls when the diffusion term depends on the control variable, even for the convex control constraint case.

The purpose of this paper is to establish some pointwise second-order necessary conditions for stochastic optimal controls in the general case.
As we shall see, both our results and the technique to prove them are quite different from that in \cite{Bonnans12, Tang10}. To see this, let us recall that, (even in the deterministic setting) in order to derive pointwise necessary conditions for optimal controls, one needs to establish first some suitable integral-type necessary conditions. In the present case,
the solution to the first-order variational equation enters into the second-order integral-type condition. Since the diffusion term contains the control variable,  there exists a term of order $\frac{3}{2}$ with respect to $\varepsilon$ (as $\varepsilon\to 0^+$) in the integral-type condition when the optimal control is perturbed by a measurable set with measure $\varepsilon$. Consequently,
the Lebesgue differentiation theorem cannot be used directly to derive the desired pointwise-type condition from the integral-type one. This is the main difficulty to treat the case that the diffusion term depends on the control variable. It is also the key trouble to derive the pointwise second-order conditions even for the case of convex control constraint.  We overcome this difficulty by means of some technique from the Malliavin calculus. On the other hand, when the control region is nonconvex, in order to derive the second-order necessary condition, one needs to expand the cost functional up to order four, and hence four adjoint equations need to be introduced.  Also, it seems interesting that, the correction part of the solution to the second-order adjoint equation (i.e. $Q_2$ in (\ref{secondajointequconvex}), or $q_2$ in (\ref{secondajointequ})) appears in the pointwise second-order necessary condition. We remark that, this part appear explicitly neither in the first-order necessary condition in \cite{Peng90}, nor in the second-order necessary conditions in the previous works \cite{Bonnans12,Tang10}.

The rest of this paper is organized as follows. Section 2 is of preliminary nature, in which we present some necessary notations and concepts.  In Section 3, we establish the pointwise second-order necessary conditions for stochastic optimal controls with convex control constraints. Finally, in Section 4, we derive the pointwise second-order necessary conditions for stochastic optimal controls with possibly nonconvex control constraints.

We refer to \cite{zhangH14a,zhangH14b} for the details of the proofs of the results in this paper and other results in this context.

\section{Some notations and concepts}

In this section, we list some notations and concepts which will be used in the sequel.

Denote by $\inner{\cdot}{\cdot}$ and $|\cdot|$ respectively the inner product and norm in $\mrn$ or $\mrm$, which can be identified from the contexts.
Let $\varphi: [0,T]\times \mrn\times \mrm\to \mrd$ ($d\in \mn$) be a map. If the map $(x,u)\mapsto \varphi(t,x,u)$ is twice differentiable for any $t\in [0,T]$, we denote by $\varphi_{(x,u)^2}(t,x,u)$ the Hessian of $\varphi$ (with respect to $(x,u)$) at $(t,x,u)$.

Denote by $\mb(\mx)$ the Borel $\sigma$-field of a metric space $\mx$.
For any $\alpha,\beta\in [1,+\infty)$, denote by $L_{\mf_{T}}^{\beta}(\Omega; \mrn)$ the space of $\mf_{T}$ measurable random variables $\xi$ such that $\me~|\xi|^{\beta}<+\infty$, by $L^{\beta}(\Omega\times[0,T]; \mrn)$ the space of $\mf\otimes\mb([0,T])$-measurable processes $\varphi$ such that $\|\varphi\|_{\beta}:=\big[\me\int_{0}^{T}|\varphi(t)|^{\beta}dt
\big]^{\frac{1}{\beta}} <+\infty$,
by $L_{\mmf}^{\beta}(\Omega; L^{\alpha}(0,T; \mrn))$ the space of $\mf\otimes\mb([0,T])$-measurable, $\mmf$-adapted processes $\varphi$ such that $\|\varphi\|_{\alpha,\beta}:=\big[\me~\big(\int_{0}^{T}|\varphi(t)|^{\alpha}dt\big)
^{\frac{\beta}{\alpha}}\big]^{\frac{1}{\beta}} <+\infty$,
by $L_{\mmf}^{\beta}(\Omega; C([0,T]; \mrn))$ the space of $\mf\otimes\mb([0,T])$-measurable, $\mmf$-adapted  continuous processes $\varphi$  such that $\|\varphi\|_{\infty,\beta}:=
\big[\me~\big(\sup_{t\in[0,T]}|\varphi(t)|^{\beta}\big)\big]^{\frac{1}{\beta}} <+\infty$, and by $L^{\infty}(\Omega\times[0,T]; \mrn)$ the space of $\mf\otimes\mb([0,T])$-measurable processes $\varphi$ such that $\|\varphi\|_{\infty}:=\esssup_{(\omega,t)\in \Omega\times[0,T]}|\varphi(\omega,t)| <+\infty $.

Also, let us recall some concepts from the Malliavin calculus. We refer to \cite{Nualart06} for a detailed discussion about this topic.
Denote by $\md^{1,2}(\mrn)$ the subspace of $L_{\mf_{T}}^{2}(\Omega; \mrn)$ whose elements are Malliavin differentiable, by $\dd_{\cdot}\xi$ the Malliavin derivative of a random variable $\xi\in \md^{1,2}(\mrn)$,
and by $\ml_{\mmf}^{1,2}(\mrn)$ the space of $\mmf$-adapted  processes $\varphi\in L_{\mmf}^{2}(\Omega;L^{2}(0,T;\mrn))$ such that
\begin{enumerate}[{\rm (i)}]
  \item \quad\quad $\varphi(t,\cdot)\in \md^{1,2}(\mrn)$, for a.e. $t\in[0,T]$;
  \item \quad\quad $(\omega,t, s)\to D_{s}\varphi(\omega,t)$ admits an $\mf\otimes\mb([0,T]\times[0,T])$-measurable version; and
  \item \quad\quad $\displaystyle|||\varphi|||_{1,2}:=\me\int_{0}^{T}|\varphi(t)|^2dt
      +\me\int_{0}^{T}\int_{0}^{T}|D_{s}\varphi(t)|^2dsdt<+\infty.$
\end{enumerate}

Further, we denote by  $\ml_{2,\mmf}^{1,2}(\mrn)$ the subspace of the stochastic processes in $\ml_{\mmf}^{1,2}(\mrn)$ whose Malliavin derivatives have suitable continuity on some neighbourhood of $\{(t,t)\big|\ t\in [0,T]\}$, i.e.,
\begin{eqnarray*}
\ml_{2,\mmf}^{1,2}(\mrn)\!\!\!\!\!&:=\!\!\!\!\!&\Big\{\varphi(\cdot)\in\ml_{\mmf}^{1,2}(\mrn)\ \Big|\ \exists\ \nabla\varphi(\cdot)\in L^2(\Omega\times[0,T];\mrn)\ \mbox{such that the functions }\\
& &\qquad f_{\ephs}:[0,T]\to [0,\infty] \mbox{ defined by } f_{\ephs}(s):=\sup_{s<t<(s+\varepsilon)\wedge T}
\me~\big|\dd_{s}\varphi(t)-\nabla\varphi(s)\big|^2,\\
& &\qquad
\mbox{for any } s\in [0,T]\mbox{ and } \ephs>0, \mbox{ are integrable, and } \lim_{\varepsilon\to 0}\int_{0}^{T}f_{\ephs}(s)ds=0\Big\}.
\end{eqnarray*}

Examples of such processes can be found in \cite{Nualart06}.

Now, let us introduce the concept of stochastic singular control. As its deterministic counterpart, a stochastic singular control is an admissible control which satisfies the first-order necessary condition trivially. Thus, before defining the stochastic singular control, let us first recall the first-order necessary condition for stochastic optimal controls established in \cite{Peng90}. Suppose that $(\bar{x}(\cdot),\bar{u}(\cdot))$ is an optimal pair.  For $\varphi=b,\sigma$ and $f$, denote
\begin{center}
\setlength{\tabcolsep}{0.5pt}
\begin{tabular*}{13cm}{@{\extracolsep{\fill}}lllr}
& $\varphi_{x}(t)=\varphi_{x}(t,\bar{x}(t),\bar{u}(t))$,
& $\varphi_{u}(t)=\varphi_{u}(t,\bar{x}(t),\bar{u}(t))$,
& $\varphi_{xx}(t)=\varphi_{xx}(t,\bar{x}(t),\bar{u}(t))$,\\
& $\varphi_{xu}(t)=\varphi_{xu}(t,\bar{x}(t),\bar{u}(t))$,
& $\varphi_{uu}(t)=\varphi_{uu}(t,\bar{x}(t),\bar{u}(t))$.
& ~~\\
\end{tabular*}
\end{center}
Define a Hamiltonian $H:\ [0,T]\times\mrn\times U\times\mrn\times\mrn\to \mr$ by
\begin{equation}\label{Hamiltonianconvex}
H(t,x,u, P,Q)
=\inner{P}{b(t,x,u)}+\inner{Q}{\sigma(t,x,u)}-f(t,x,u),
\end{equation}
for any $(t,x,u,P,Q)\in [0,T]\times\mrn\times U\times\mrn\times\mrn$.
Let $(P_{1}(\cdot),Q_{1}(\cdot))$ and $(P_{2}(\cdot),Q_{2}(\cdot))$ solve respectively the following two adjoint equations
\begin{equation}\label{firstajointequconvex}
 \left\{
\begin{array}{l}
dP_{1}(t)=-\Big[b_{x}(t)^{\top}P_{1}(t)+\sigma_{x}(t)^{\top}Q_{1}(t)
          -f_{x}(t)\Big]dt+Q_{1}(t)dW(t), \  t\in[0,T], \\
P_{1}(T)=-h_{x}(\bar{x}(T))
\end{array}\right.
\end{equation}
and
\begin{equation}\label{secondajointequconvex}
\quad\left\{
\begin{array}{l}
dP_{2}(t)=-\Big[b_{x}(t)^{\top}P_{2}(t)+P_{2}(t)b_{x}(t) +\sigma_{x}(t)^{\top}P_{2}(t)\sigma_{x}(t) +\sigma_{x}(t)^{\top}Q_{2}(t)\\
\qquad\qquad \qquad\qquad \qquad\quad\quad
+Q_{2}(t)\sigma_{x}(t)+H_{xx}(t)\Big]dt+Q_{2}(t)dW(t),\  t\in[0,T], \\
P_{2}(T)=-h_{xx}(\bar{x}(T)),
\end{array}\right.
\end{equation}
where
$H_{xx}(t)=H_{xx}(t,\bar{x}(t),\bar{u}(t), P_{1}(t),Q_{1}(t))$.
Define another function $\mh: [0,T]\times\mrn\times U\to\mr$ by
\begin{eqnarray}\label{mathbb H}
& &\mh(t,x,u)=H(t,x,u, P_{1}(t),Q_{1}(t))-H(t,x,\bar{u}(t), P_{1}(t),Q_{1}(t))\nonumber\\
& &\qquad\qquad\quad +\frac{1}{2}\inner{P_{2}(t)\big(\sigma(t,x,u)-\sigma(t,x,\bar{u}(t))\big)}
{\sigma(t,x,u)-\sigma(t,x,\bar{u}(t))},\nonumber\\
& &\qquad\qquad\qquad\qquad\qquad\qquad\qquad\qquad\qquad\qquad\qquad\qquad\quad
(t,x,u)\in [0,T]\times\mrn\times U.
\end{eqnarray}
It was shown in \cite{Peng90} that the optimal pair $(\bar{x}(\cdot),\bar{u}(\cdot))$ satisfies
\begin{equation}\label{peng's maximum principle}
\mh(t,\bar{x}(t),v)\le 0,\quad \forall\ v\in U,
\ a.e.\ (\omega,t)\in \Omega\times [0,T].
\end{equation}
That is, the function $v\mapsto \mh(t,\bar{x}(t),v)$ attends its maximum at $\bar{u}(t)$ for a.e. $(\omega,t)\in \Omega\times [0,T]$. By the first- and second-order necessary conditions in classical optimization theory, when $U$ is convex and $b,\ \sigma$ and $f$ are sufficiently smooth, $\bar{u}(\cdot)$ satisfies
\begin{eqnarray}\label{first neces condi for maximum principle}
\inner{\mh_{u}(t,\bar{x}(t),\bar{u}(t))}{v-\bar{u}(t)}
&=&\inner{H_{u}(t,\bar{x}(t), \bar{u}(t) ,P_{1}(t),Q_{1}(t))}{v-\bar{u}(t)}\nonumber\\
&\le& 0,\quad
\forall\ v\in U, \ a.e.\  (\omega,t)\in \Omega\times [0,T].
\end{eqnarray}
Moreover, if
$H_{u}(t,\bar{x}(t), \bar{u}(t) ,P_{1}(t),Q_{1}(t))= 0,\ a.e.\ (\omega,t)\in \Omega\times [0,T],$
then
\begin{eqnarray}\label{second neces condi for maximum principle}
& &\inner{\mh_{uu}(t,\bar{x}(t),\bar{u}(t))(v-\bar{u}(t))}{v-\bar{u}(t)}\nonumber\\
& &=\inner{H_{uu}(t,\bar{x}(t), \bar{u}(t) ,P_{1}(t),Q_{1}(t))(v-\bar{u}(t))}{v-\bar{u}(t)}\nonumber\\
& &\ \ \ + \inner{\sigma_{u}(t,\bar{x}(t),\bar{u}(t))^{\top}
P_{2}(t)\sigma_{u}(t,\bar{x}(t),\bar{u}(t))(v-\bar{u}(t)) }{v-\bar{u}(t)}\nonumber\\
& &\le 0,\quad
\forall\ v\in U, \ a.e.\ (\omega,t)\in \Omega\times [0,T].
\end{eqnarray}

According to (\ref{first neces condi for maximum principle})--(\ref{second neces condi for maximum principle}) and (\ref{peng's maximum principle}), we give below the definition of stochastic singular control in the classical sense and that in the sense of Pontryagin-type maximum principle, respectively.

\begin{definition}
Let $\tilde{u}(\cdot)\in \mmu_{ad}$, $\tilde{x}(\cdot)$ be the state with respect to $\tilde{u}(\cdot)$, and $(\tilde{P}_{1}(\cdot),\tilde{Q}_{1}(\cdot))$ and $(\tilde{P}_{2}(\cdot),\tilde{Q}_{2}(\cdot))$ be the adjoint processes given respectively by (\ref{firstajointequconvex}) and (\ref{secondajointequconvex}) with $(\bar x(\cdot),\bar u(\cdot)
)$ replaced by $(\tilde{x}(\cdot),\tilde{u}(\cdot))$.
\begin{enumerate}[{\rm (i)}]
\item The control $\tilde{u}(\cdot)$ is said to be singular in the classical sense if, for a.e. $(\omega,t)\in \Omega\times [0,T]$,
      \begin{equation}\label{singularcontrol calssical}
      \left\{
        \begin{array}{l}
         H_{u}(t,\tilde{x}(t), \tilde{u}(t) ,\tilde{P}_{1}(t),\tilde{Q}_{1}(t))=0, \\
         H_{uu}(t,\tilde{x}(t), \tilde{u}(t) ,\tilde{P}_{1}(t),\tilde{Q}_{1}(t))+\sigma_{u}(t,\tilde{x}(t),\tilde{u}(t))^{\top}
        \tilde{P}_{2}(t)\sigma_{u}(t,\tilde{x}(t),\tilde{u}(t))=0;
        \end{array}\right.
\end{equation}

\item The control $\tilde{u}(\cdot)$ is said to be singular in the sense of Pontryagin-type maximum principle on a control region $V$ if, $V$ is a nonempty subset of $U$ and for a.e. $(\omega,t)\in \Omega\times [0,T]$,
    \begin{eqnarray}\label{singularity for maximum principle}
      & &H(t,\tilde{x}(t), v, \tilde{P}_{1}(t), \tilde{Q}_{1}(t))
      -H(t,\tilde{x}(t), \tilde{u}(t), \tilde{P}_{1}(t),\tilde{Q}_{1}(t))\nonumber\\
      & &+\frac{1}{2}\inner{\tilde{P}_{2}(t)\big(\sigma(t,\tilde{x}(t),v)                         -\sigma(t,\tilde{x}(t),\tilde{u}(t))\big)}                         {\sigma(t,\tilde{x}(t),v)-\sigma(t,\tilde{x}(t),\tilde{u}(t))}
      =0.
\end{eqnarray}
\end{enumerate}

\end{definition}

\begin{remark}
The above two concepts of stochastic singular controls are natural extensions of their deterministic counterpart (see \cite{Gabasov72}). Also, when the set $U$ is open and $V=U$, any admissible control satisfying
(\ref{singularity for maximum principle}) must satisfy
(\ref{singularcontrol calssical}). In this case, every admissible control which is singular in the sense of the Pontryagin-type maximum principle is also singular in the classical sense.
\end{remark}

\section{Pointwise second-order necessary conditions, the convex control constraint case}

In this section, we consider the pointwise second-order necessary condition for stochastic optimal controls with convex control constraints. Here, the optimal controls are assumed to be singular in the classical sense.

Similar to \cite{Bonnans12}, we assume that
\begin{enumerate}
  \item [{\em (C1)}] The control region $U $ is nonempty, bounded and convex.
  \item [{\em (C2)}] The maps $b$, $\sigma$, $f$ and $h$ satisfy the following:
  \begin{enumerate}[{\rm (i)}]
       \item $b$, $\sigma$ and $f$ are
             $\mb([0, T]\times\mrn\times U)$-measurable, $h$ is $\mb(\mrn)$-measurable.
       \item For a.e. $t\in [0, T]$, the map
             $(x,u)\mapsto (b(t, x, u), \sigma(t, x, u))$
             has continuous bounded partial derivatives up to order three. And, there exists a constant $L > 0$ such that $|\varphi(t,0, u)|\le L$, for $\varphi = b,\ \sigma$, a.e. $t\in [0,T]$ and any $u\in U$.

       \item For a.e. $t\in [0, T]$, the map
             $(x,u)\mapsto (f(t, x, u), h(x))$
             has continuous partial derivatives up to order three. And, there exists a constant $L > 0$ such that for a.e.
             $t\in [0,T]$ and any $x,\ \tilde{x}\in\mrn$, $u,\ \tilde{u}\in U$, it holds that
             \begin{equation}\nonumber
             \left\{
             \begin{array}{l}
              |f(t,x,u)|\le L(1+|x|^{2}+|u|^{2}),\\
              |f_{x}(t,x,u)|+|f_{u}(t,x,u)|\le L(1+|x|+|u|),\\
              |f_{xx}(t,x,u)|+|f_{xu}(t,x,u)|+|f_{uu}(t,x,u)|\le L,\\
              |f_{(x,u)^2}(t,x,u)-f_{(x,u)^2}(t,\tilde{x},\tilde{u})|\le L(|x-\tilde{x}|+|u-\tilde{u}|),\\
              |h(x)|\le L(1+|x|^{2}),\ |h_{x}(x)| \le L(1+|x|),\\
              |h_{xx}(x)| \le L, \ |h_{xx}(x)-h_{xx}(\tilde{x})|\le L|x-\tilde{x}|.
             \end{array}\right.
             \end{equation}
               \end{enumerate}
\end{enumerate}

Since the control region is assumed to be convex, the convex variation  can be used as a perturbation of the optimal control.
Let $(\bar{x}(\cdot),\bar{u}(\cdot))$ be an optimal pair, and $u(\cdot)\in \mathcal{U}_{ad}$ be any given admissible control.  For any  $\varepsilon\in(0,1)$, write
\begin{equation}\nonumber
v(\cdot)=u(\cdot)-\bar{u}(\cdot),\qquad u^{\varepsilon}(\cdot)=\bar{u}(\cdot)+\varepsilon v(\cdot).
\end{equation}
Clearly, $u^{\varepsilon}(\cdot)\in \mmu_{ad}$.
Denote by $x^{\varepsilon}(\cdot)$ the state with respect to the control
$u^{\varepsilon}(\cdot)$, and put $\delta x(\cdot)=x^{\varepsilon}(\cdot)-\bar{x}(\cdot)$.

We introduce the following two variational equations:
\begin{equation}\label{variational equation one}
\quad\left\{
\begin{array}{l}
dy_{1}(t)= \Big[b_{x}(t) y_{1}(t)+b_{u}(t)v(t)\Big]dt
                +\Big[\sigma_{x}(t) y_{1}(t)+ \sigma_{u}(t)v(t)\Big]dW(t),\  t\in[0,T],\qquad\quad\qquad\qquad\quad\\
y_{1}(0)=0
\end{array}\right.
\end{equation}
and
\begin{equation}\label{variational equation two}
\quad\left\{
\begin{array}{l}
dy_{2}(t)= \Big[b_{x}(t)
y_{2}(t)+y_{1}(t)^{\top}b_{xx}(t)y_{1}(t)+2v(t)^{\top}b_{xu}(t)y_{1}(t)
+v(t)^{\top}b_{uu}(t)v(t)\Big]dt\\ [+0.5em]
\qquad\quad
+\Big[\sigma_{x}(t)y_{2}(t)+y_{1}(t)^{\top}\sigma_{xx}(t)y_{1}(t)
+2v(t)^{\top}\sigma_{xu}(t)y_{1}(t)
+v(t)^{\top}\sigma_{uu}(t)v(t)\Big]dW(t),\ t\in[0,T],\\
y_{2}(0)=0.
\end{array}\right.
\end{equation}
The solutions $y_{1}(\cdot)$ and $y_{2}(\cdot)$ are the first- and the second-order linear approximations of the state variation $\delta x(\cdot)$, respectively. By the duality between the variational equations (\ref{variational equation one})--(\ref{variational equation two}) and the adjoint equations (\ref{firstajointequconvex})--(\ref{secondajointequconvex}), the following integral-type second-order necessary condition immediately follows.

\begin{theorem}
Let  (C1)--(C2) hold. If $\bar{u}(\cdot)$ is a singular optimal control in the classical sense, then
\begin{equation}\label{integraltype 2ordercondition}
\me\int^{T}_{0}\inner{\ms(t)y_{1}(t)}{v(t)}dt\le 0,  \quad \forall \ v(\cdot)=u(\cdot)-\bar{u}(\cdot),\ u(\cdot)\in \mmu_{ad},
\end{equation}
where the process $\ms(\cdot)$ is defined by
\begin{equation}\label{S(t)}
\ms(t)
:= H_{xu}(t)+b_{u}(t)^{\top}P_{2}(t)
 +\sigma_{u}(t)^{\top}Q_{2}(t)+\sigma_{u}(t)^{\top}
P_{2}(t)\sigma_{x}(t),\quad t\in[0,T],
\end{equation}
and $H_{xu}(t)=H_{xu}(t,\bar{x}(t),\bar{u}(t),P_{1}(t),Q_{1}(t))$.
\end{theorem}

Next, we will rewrite the integral-type necessary condition (\ref{integraltype 2ordercondition}) as a pointwise one.

Fix $\t\in [0,T)$ arbitrarily. Define
$$u(t)=\left\{
\begin{array}{l}
v, \qquad\qquad t\in E_{\theta},\\
\bar{u}(t), \qquad \quad t\in [0,T] \setminus E_{\theta},\\
\end{array}\right.
$$
where $v\in U$, $E_{\theta}=[\t,\t+\theta)$, $\theta>0$, $\t+\theta\le T$. Denote by $\chi_{E_{\theta}}(\cdot)$ the characteristic function of the set $E_{\theta}$, and by $\Phi(\cdot)$ the solution to the following matrix-value stochastic differential equation:
\begin{equation*}
\left\{
\begin{array}{l}
d\Phi(t)= b_{x}(t)\Phi(t)dt+\sigma_{x}(t)\Phi(t)dW(t),
\qquad \ \ \ t\in[0,T], \qquad \\
\Phi(0)=I,
\end{array}\right.
\end{equation*}
where $I$ is the identity matrix in $\mr^{n\times n}$. By \cite[Theorem 1.6.14]{Yong99} (at page 47), the solution $y_{1}(\cdot)$ to the equation (\ref{variational equation one}) with respect to $v(\cdot)=u(\cdot)-\bar{u}(\cdot)=(v-\bar{u}(\cdot))\chi_{E_{\theta}}(\cdot)$ enjoys the following explicit representation:
\begin{eqnarray}\label{y1(t)for v-u(t)}
 y_{1}(t)\!\!\!&=\!\!\!&\Phi(t)\int_{0}^{t}\Phi(s)^{-1}
\big(b_{u}(s)(v-\bar{u}(s))-\sigma_{x}(s)\sigma_{u}(s)(v-\bar{u}(s))
\big)\chi_{E_{\theta}}(s)ds\nonumber\\
& &+\Phi(t)\int_{0}^{t}
\Phi(s)^{-1}\sigma_{u}(s)(v-\bar{u}(s))\chi_{E_{\theta}}(s)dW(s), \quad t\in[0,T].
\end{eqnarray}

Substituting (\ref{y1(t)for v-u(t)}) into (\ref{integraltype 2ordercondition}), there will appear the following term
\begin{equation}\label{3 per 2 order term}
\me\int_{\t}^{\t+\theta}
\Big\langle\ms(t)\Phi(t)\int_{\t}^{t}
\Phi(s)^{-1}\sigma_{u}(s)
(v-\bar{u}(s))dW(s),
v-\bar{u}(t)\Big\rangle dt,
\end{equation}
which is not an infinitesimal of order two but only that of order $\frac{3}{2}$ with respect to $\theta$ (as $\theta\to 0^+$).
Therefore, the Lebesgue differentiation theorem cannot be used directly to derive the pointwise second-order necessary condition (from the integral-type second order necessary condition \eqref{integraltype 2ordercondition}). To overcome this difficulty, we need to introduce the following regularity assumption:

\begin{enumerate}
\item [{\em (C3)}]
$$\bar{u}(\cdot)\in\ml_{2,\mmf}^{1,2}(\mrm),\mbox{ and } \ms(\cdot)\in \ml_{2,\mmf}^{1,2}(\mr^{m\times n})\cap L^{\infty}(\Omega\times[0,T]; \mr^{m\times n}).$$
\end{enumerate}

By the boundness of $U$ and the assumption {\em (C3)}, it follows that
$\ms(\cdot)^{\top}(v-\bar{u}(\cdot))\in\ml^{1,2}_{\mmf}(\mrn)\cap L^{\infty}(\Omega\times[0,T];\mr^{n})$, for any $v\in U$.
By the Clark-Ocone formula, we have
\begin{equation}\label{expu(t)S(t)}
\ms(t)^{\top}(v-\bar{u}(t))
=\me~\Big[\ms(t)^{\top}(v-\bar{u}(t))\Big]
 +\int_{0}^{t}\me~\Big[\dd_{s}
\big(\ms(t)^{\top}(v-\bar{u}(t))\big)\;\Big|\;\mf_{s}\Big]dW(s).
\end{equation}
Substituting (\ref{expu(t)S(t)}) into (\ref{3 per 2 order term}), and using the properties of the It\^{o} integral and the conditional expectation, we see that
\begin{eqnarray*}
& &\me\int_{\t}^{\t+\theta}
\Big\langle\ms(t)\Phi(t)\int_{\t}^{t}
\Phi(s)^{-1}\sigma_{u}(s)
(v-\bar{u}(s))dW(s),
v-\bar{u}(t)\Big\rangle dt\\
& &=\me\int_{\t}^{\t+\theta}\int_{\t}^{t}
\Big\langle\Phi(\t)\Phi(s)^{-1}\sigma_{u}(s)
(v-\bar{u}(s)),
\dd_{s}\big(\ms(t)^{\top}(v-\bar{u}(t))\big)
\Big\rangle dsdt\\
& &\ \ \ +\me\int_{\t}^{\t+\theta}
\int_{\t}^{t}\Big\langle\sigma_{x}(s)\sigma_{u}(s)(v-\bar{u}(s)),
\ms(\t)^{\top}\big((v-\bar{u}(\t)\big)\Big\rangle dsdt
+o(\theta^2), \quad(\mbox{as }\theta\to 0^+).
\end{eqnarray*}
In this way, the Lebesgue differentiation theorem can be used and then we obtain the following pointwise secend-order necessary condition.
\begin{theorem}\quad
\label{2orderconditionth}
Let (C1)--(C3) hold. If $\bar{u}(\cdot)$ is a singular optimal control in the classical sense, then for a.e. $\t\in [0,T]$, it holds that
\begin{eqnarray}\label{2orderconditionconvex}
& &\inner{\ms(\t)b_{u}(\t)(v-\bar{u}(\t))}{v-\bar{u}(\t)}\nonumber\\
& &\quad+ \inner{\nabla \ms(\t)\sigma_{u}(\t)(v-\bar{u}(\t))}
{v-\bar{u}(\t)}- \inner{\ms(\t)\sigma_{u}(\t)(v-\bar{u}(\t))}
{\nabla\bar{u}(\t)}\\
& &
\le 0, \quad\quad \forall \ v\in U, \ a.s.\nonumber
\end{eqnarray}
\end{theorem}

\section{Pointwise second-order necessary conditions, the general case}
In this section, we discuss the pointwise second-order necessary condition for stochastic optimal controls with possibly nonconvex control constraints. Here, the optimal controls are assumed to be singular in the sense of Pontryagin-type maximum principle. Unlike the convex control constraint case, the cost functional needs to be expanded up to order four, and therefore four variational equations and four adjoint equations need to be introduced. To avoid introducing high order tensors, we only consider the 1-dimensional case here, i.e., $m=n=1$, and hence both the control and the state are assumed to take values in $\mr$.

In this section, we assume that
\begin{enumerate}
    \item [{\em (C4)}] The control region $U \subset \mr$ is nonempty and bounded.
    \item [{\em (C5)}] The maps $b,\ \sigma,$  $f$, and $h$ satisfy the following:
       \begin{enumerate}[{\rm (i)}]
              \item $b$, $\sigma$ and $f$ are $\mb([0; T]\times\mr\times U)$-measurable, $h$ is $\mb(\mr)$-measurable.
              \item For a.e. $(t, u)\in [0; T] \times U$, the map
                 $x\mapsto (b(t, x, u), \sigma(t, x, u), f(t, x, u))$
                is continuously differentiable up to order four, and
                there exist a constant $L > 0$ and a modulus of continuity $\tilde{\omega}: [0,\infty) \to [0, \infty)$ such that for $\varphi = b,\ \sigma,\ f$ it holds that, for a.e. $ t\in [0,T]$, and any $x,\ \tilde{x}\in\mr$, $u,\ \tilde{u}\in U$,
                   \begin{equation}\nonumber
                   \left\{
                    \begin{array}{l}
                      |\varphi(t,x, u)-\varphi(t,\tilde{x}, \tilde{u})|\le L|x-\tilde{x}|+ \tilde{\omega}(|u-\tilde{u}|),\\
                          |\varphi_{x}(t,x, u)-\varphi_{x}(t,\tilde{x}, \tilde{u})|
                      \le L|x-\tilde{x}|+ \tilde{\omega}(|u-\tilde{u}|),\\
                     |\varphi_{xx}(t,x, u)-\varphi_{xx}(t,\tilde{x}, \tilde{u})|
                      \le L|x-\tilde{x}|+ \tilde{\omega}(|u-\tilde{u}|),\\
                    |\varphi_{xxx}(t,x, u)-\varphi_{xxx}(t,\tilde{x}, \tilde{u})|
                    \le L|x-\tilde{x}|+ \tilde{\omega}(|u-\tilde{u}|),\\
                    |\varphi_{xxxx}(t,x, u)-\varphi_{xxxx}(t,x, \tilde{u})|\le L|x-\tilde{x}|+\tilde{\omega}(|u-\tilde{u}|).\\
                  \end{array}\right.
                    \end{equation}
                \item $h$ is continuously differentiable up to order four, and
                there exists a constant $L > 0$ such that for any $x\in \mr$,
                   \begin{eqnarray*}
                   & &|h(x)|\le L(1+|x|^{4}), \quad\ \ \  |h_{x}(x)|\le L(1+|x|^3),\quad \\
                   & &|h_{xx}(x)|\le L(1+|x|^{2}), \quad |h_{xxx}(x)|\le L(1+|x|),\quad |h_{xxxx}(x)|\le L.
                   \end{eqnarray*}
             \end{enumerate}
\end{enumerate}

Firstly, we establish a variational formulation for optimal controls. Since the control region may be nonconvex, we need to use the needle variation as a perturbation of the optimal control.

Let $(\bar{x}(\cdot),\bar{u}(\cdot))$ be an optimal pair, $v\in U$, and $E_{\varepsilon} \subset [0,T]$ be a Lebesgue measurable set satisfies $|E_{\varepsilon}|=\varepsilon$, where $|E_{\varepsilon}|$ stands for the Lebesgue measure of $E_{\varepsilon}$. Define
\begin{equation}\label{needle variation}
u^{\varepsilon}(t)=
\begin{cases}
v, &\text{$t\in E_\varepsilon$},\\
\bar{u}(t), &\text{$t\in [0,T]\setminus E_{\varepsilon}$}.
\end{cases}
\end{equation}
Obviously, $u^{\varepsilon}(\cdot)\in \mmu_{ad}$.
Let $x^{\varepsilon}(\cdot)$ be the state with respect to the control $u^{\varepsilon}(\cdot)$. Denote $\delta x(\cdot)=x^{\varepsilon}(\cdot)-\bar{x}(\cdot)$ and, for $\varphi=b,\sigma$ and $f$, write
\begin{center}
\setlength{\tabcolsep}{0.5pt}
\begin{tabular*}{16cm}{@{\extracolsep{\fill}}lllr}
& $\delta \varphi(t)=\varphi(t,\bar{x}(t),v)
-\varphi(t,\bar{x}(t),\bar{u}(t))$,
& $\delta \varphi_{x}(t)=\varphi_{x}(t,\bar{x}(t),v)
-\varphi_{x}(t,\bar{x}(t),\bar{u}(t))$\\
& $\delta \varphi_{xx}(t)=\varphi_{xx}(t,\bar{x}(t),v)
-\varphi_{xx}(t,\bar{x}(t),\bar{u}(t))$,
& $\delta\varphi_{xxx}(t)=\varphi_{xxx}(t,\bar{x}(t),v)
-\varphi_{xxx}(t,\bar{x}(t),\bar{u}(t))$.
\end{tabular*}
\end{center}
\vspace{2mm}

We introduce the following four variational equations:
\begin{equation}\label{firstvariequ}
\left\{
\begin{array}{l}
dy_{1}^{\varepsilon}(t)= b_{x}(t) y_{1}^{\varepsilon}(t)dt
+\Big[\sigma_{x}(t) y_{1}^{\varepsilon}(t)+ \delta\sigma(t)\chi_{E_{\varepsilon}}(t)\Big]dW(t),\quad t\in [0,T], \qquad \\
y_{1}^{\varepsilon}(0)=0,
\end{array}\right.
\end{equation}
\begin{equation}\label{secondvariequ}
\left\{
\begin{array}{l}
dy_{2}^{\varepsilon}(t)= \Big[b_{x}(t)
y_{2}^{\varepsilon}(t)+\frac{1}{2}b_{xx}(t)y_{1}^{\varepsilon}(t)^{2}+\delta b(t)\chi_{E_{\varepsilon}}(t)\Big]dt\\[+0.3em]
\qquad\ \  +\Big[\sigma_{x}(t) y_{2}^{\varepsilon}(t)
+\frac{1}{2}\sigma_{xx}(t)y_{1}^{\varepsilon}(t)^{2}
+\delta\sigma_{x}(t)y_{1}^{\varepsilon}(t)\chi_{E_{\varepsilon}}(t)\Big]dW(t),\ t\in [0,T], \\
y_{2}(0)=0,
\end{array}\right.
\end{equation}
\begin{equation}\label{thirdvariequ}
\left\{
\begin{array}{l}
dy_{3}^{\varepsilon}(t)= \Big[b_{x}(t)y_{3}^{\varepsilon}(t)
+\frac{1}{2}b_{xx}(t)\big(2 y_{1}^{\varepsilon}(t) y_{2}^{\varepsilon}(t)+ y_{2}^{\varepsilon}(t)^{2}\big)\\[+0.2em]
\qquad\qquad\qquad+\frac{1}{6}b_{xxx}(t) y_{1}^{\varepsilon}(t)^{3}
+\delta b_{x}(t) y_{1}^{\varepsilon}(t) \chi_{E_{\varepsilon}}(t)\Big]dt\\[+0.2em]
\qquad \qquad +\Big[\sigma_{x}(t) y_{3}^{\varepsilon}(t)
+ \frac{1}{2}\sigma_{xx}(t)\big(2 y_{1}^{\varepsilon}(t) y_{2}^{\varepsilon}(t)+ y_{2}^{\varepsilon}(t)^{2}\big)+\frac{1}{6}\sigma_{xxx}(t) y_{1}^{\varepsilon}(t)^{3}\ \ \ \ \\[+0.2em]
\qquad \qquad+\delta \sigma_{x}(t) y_{2}^{\varepsilon}(t) \chi_{E_{\varepsilon}}(t) +\frac{1}{2}\delta\sigma_{xx}(t)y_{1}^{\varepsilon}(t)^{2}
\chi_{E_{\varepsilon}}(t)\Big]dW(t),\quad t\in [0,T],\\[+0.2em]
y_{3}^{\varepsilon}(0)=0,
\end{array}\right.
\end{equation}
and
\begin{equation}\label{forthvariequ}
\left\{
\begin{array}{l}
dy_{4}^{\varepsilon}(t)= \Big[b_{x}(t) y_{4}^{\varepsilon}(t)
+\frac{1}{2}b_{xx}(t)\big(2 y_{1}^{\varepsilon}(t) y_{3}^{\varepsilon}(t)+ 2 y_{2}^{\varepsilon}(t) y_{3}^{\varepsilon}(t)+ y_{3}^{\varepsilon}(t)^{2}\big)\\[+0.2em]
\qquad\qquad +\frac{1}{6}b_{xxx}(t)\big (3y_{1}^{\varepsilon}(t)^{2}
y_{2}^{\varepsilon}(t)+3y_{1}^{\varepsilon}(t) y_{2}^{\varepsilon}(t)^{2}
+y_{2}^{\varepsilon}(t)^{3}\big)\\[+0.2em]
\qquad\qquad+\frac{1}{24}b_{xxxx}(t) y_{1}^{\varepsilon}(t)^{4}
+\delta b_{x}(t) y_{2}^{\varepsilon}(t) \chi_{E_{\varepsilon}}(t)+\frac{1}{2}\delta b_{xx}(t) y_{1}^{\varepsilon}(t)^{2}\chi_{E_{\varepsilon}}(t)\Big]dt\\[+0.2em]
\qquad\qquad +\Big[\sigma_{x}(t) y_{4}^{\varepsilon}(t)
+ \frac{1}{2}\sigma_{xx}(t)\big(2 y_{1}^{\varepsilon}(t) y_{3}^{\varepsilon}(t)+ 2 y_{2}^{\varepsilon}(t) y_{3}^{\varepsilon}(t)+ y_{3}^{\varepsilon}(t)^{2}\big)\\[+0.2em]
\qquad \qquad+\frac{1}{6}\sigma_{xxx}(t) \big(3y_{1}^{\varepsilon}(t)^{2} y_{2}^{\varepsilon}(t)+3y_{1}^{\varepsilon}(t) y_{2}^{\varepsilon}(t)^{2}
+y_{2}^{\varepsilon}(t)^{3}\big)\\[+0.2em]
\qquad \qquad +\frac{1}{24}\sigma_{xxxx}(t) y_{1}^{\varepsilon}(t)^{4}
+\delta \sigma_{x}(t) y_{3}^{\varepsilon}(t) \chi_{E_{\varepsilon}}(t) +\frac{1}{2}\delta\sigma_{xx}(t)(2 y_{1}^{\varepsilon}(t) y_{2}^{\varepsilon}(t) \\[+0.2em]
\qquad\qquad
+ y_{2}^{\varepsilon}(t)^{2})\chi_{E_{\varepsilon}}(t)
+\frac{1}{6}\delta\sigma_{xxx}(t)
y_{1}^{\varepsilon}(t)^{3}\chi_{E_{\varepsilon}}(t)\Big]dW(t), \quad t\in [0,T],\\[+0.2em]
y_{4}^{\varepsilon}(0)=0.
\end{array}\right.
\end{equation}

Corresponding to the variational equations (\ref{firstvariequ})--(\ref{forthvariequ}), we introduce the following four adjoint equations:
\begin{equation}\label{firstajointequ}
\left\{
\begin{array}{l}
dp_{1}(t)=-\Big[b_{x}(t)p_{1}(t)+\sigma_{x}(t)q_{1}(t)-f_{x}(t)\Big]dt+q_{1}(t)dW(t),
 \  t\in[0,T],\qquad\qquad\qquad\qquad \\
p_{1}(T)=-h_{x}(\bar{x}(T)),
\end{array}\right.
\end{equation}
\begin{equation}\label{secondajointequ}
\left\{
\begin{array}{l}
dp_{2}(t)=-\Big[2b_{x}(t)p_{2}(t) +\sigma_{x}(t)^{2}p_{2}(t)
+2\sigma_{x}(t)q_{2}(t)
+H_{xx}(t)\Big]dt+q_{2}(t)dW(t),\  t\in[0,T],\  \\
p_{2}(T)=-h_{xx}(\bar{x}(T)),
\end{array}\right.
\end{equation}
\begin{equation}\label{thirdajointequ}
\left\{
\begin{array}{l}
dp_{3}(t)=-\Big[3b_{x}(t)p_{3}(t)+3\sigma_{x}^{2}(t)p_{3}(t)+3\sigma_{x}(t)q_{3}(t)
+3b_{xx}(t)p_{2}(t) +3\sigma_{xx}(t)q_{2}(t)\qquad\qquad \\
\qquad\qquad\qquad+3\sigma_{x}(t)\sigma_{xx}(t)p_{2}(t)
+H_{xxx}(t)\Big]dt+q_{3}(t)dW(t),\ t\in [0,T],  \\
p_{3}(T)=-h_{xxx}(\bar{x}(T)),
\end{array}\right.
\end{equation}
and
\begin{equation}\label{forthajointequ}
\left\{
\begin{array}{l}
dp_{4}(t)=-\Big[4b_{x}(t)p_{4}(t)+6\sigma_{x}^{2}(t)p_{4}(t)
+4\sigma_{x}(t)q_{4}(t)+6b_{xx}(t)p_{3}(t) +6\sigma_{xx}(t)q_{3}(t)\\
\qquad\qquad\qquad+12\sigma_{x}(t)\sigma_{xx}(t)p_{3}(t)
+4b_{xxx}(t)p_{2}(t)+4\sigma_{x}(t)\sigma_{xxx}(t)p_{2}(t) +3\sigma_{xx}^{2}(t)p_{2}(t)\quad\ \ \\
\qquad\qquad\qquad+4\sigma_{xxx}(t)q_{2}(t)+H_{xxxx}(t)\Big]dt
+q_{4}(t)dW(t),\ t\in [0,T],\quad\ \\
p_{4}(T)=-h_{xxxx}(\bar{x}(T)),
\end{array}\right.
\end{equation}
where the Hamiltonian $H$ is defined by (\ref{Hamiltonianconvex}) (with $n=1$),
\begin{center}
\setlength{\tabcolsep}{0.5pt}
\begin{tabular*}{16cm}{@{\extracolsep{\fill}}lllr}
& $H_{xx}(t)=H_{xx}(t, \bar{x}(t),\bar{u}(t), p_{1}(t),q_{1}(t))$,
& $H_{xxx}(t)=H_{xxx}(t, \bar{x}(t),\bar{u}(t), p_{1}(t),q_{1}(t))$,\\
& $H_{xxxx}(t)=H_{xxxx}(t, \bar{x}(t),\bar{u}(t), p_{1}(t),q_{1}(t))$.
& ~~
\end{tabular*}
\end{center}
\vspace{2mm}

In addition, we define the following two functions $\mS,\ \mT: [0,T]\times\mr\times U\times\mr\times\mr\to \mr$ by
\begin{eqnarray*}
\mS(t,x,u, p_{2},q_{2})=p_{2}b(t,x,u)+q_{2}\sigma(t,x,u), \qquad \mT(t,x,u, p_{3},q_{3})=p_{3}b(t,x,u)+q_{3}\sigma(t,x,u),
 \end{eqnarray*}
for $(t,x,u,p_{3},q_{3})\in [0,T]\times\mr\times U\times\mr\times\mr$, and the functions $\mh,\ \mss,\ \mt: [0,T]\times\mr\times U\to \mr$ by
\begin{eqnarray*}
\mh(t,x,u)&=&H(t,x,u, p_{1}(t),q_{1}(t))-H(t,x,\bar{u}(t), p_{1}(t),q_{1}(t))
+\frac{1}{2}p_{2}(t)\big(\sigma(t,x,u)-\sigma(t,x,\bar{u}(t))\big)^{2},\\
\mss(t,x,u)&=&\mh_{x}(t,x,u)+\mS(t,x,u, p_{2}(t),q_{2}(t))
-\mS(t,x,\bar{u}(t), p_{2}(t),q_{2}(t))\\
& &+p_{2}\sigma_{x}(t,x,\bar{u}(t))\big(\sigma(t,x,u)
-\sigma(t,x,\bar{u}(t))\big)
+\frac{1}{2}p_{3}\big(\sigma(t,x,u)-\sigma(t,x,\bar{u}(t))\big)^{2},\\
\mt(t,x,u)&=&\mss_{x}(t,x,u)+\mS_{x}(t,x,u, p_{2}(t),q_{2}(t))
-\mS_{x}(t,x,\bar{u}(t), p_{2}(t),q_{2}(t))\\
& &+\mT(t,x,u, p_{3}(t),q_{3}(t))
   -\mT(t,x,\bar{u}(t), p_{3}(t),q_{3}(t))\\
& &+p_{2}\sigma_{x}(t,x,\bar{u}(t))\big(\sigma_{x}(t,x,u)
-\sigma_{x}(t,x,\bar{u}(t))\big)\\
& &+p_{3}\big(\sigma(t,x,u)
-\sigma(t,x,\bar{u}(t))\big)\big(\sigma_{x}(t,x,u)
-\sigma_{x}(t,x,\bar{u}(t))\big)\\
& &+2p_{3}\sigma_{x}(t,x,\bar{u}(t))\big(\sigma(t,x,u)
-\sigma(t,x,\bar{u}(t))\big)
+\frac{1}{2}p_{4}\big(\sigma(t,x,u(t))-\sigma(t,x,\bar{u}(t))\big)^{2},
\end{eqnarray*}
for $(t,x,u)\in [0,T]\times\mr\times U$.

Using the duality between the variational equations (\ref{firstvariequ})--(\ref{forthvariequ}) and the adjoint equations (\ref{firstajointequ})--(\ref{forthajointequ}), we obtain the following variational formulation for the optimal control $\bar{u}(\cdot)$.

\begin{proposition}
Let (C4)--(C5) hold.  Then the following variational equality holds:
\begin{equation}\label{shorttaylor}
\begin{split}
 J(u^{\varepsilon})-J(\bar{u})
=&-\me\int_{0}^{T}\Big[
\mh(t,\bar{x}(t),v)
+\ms(t,\bar{x}(t),v)(y_{1}^{\varepsilon}(t)+y_{2}^{\varepsilon}(t))\\
&\qquad\qquad  +\frac{1}{2}\mt(t,\bar{x}(t),v)y_{1}^{\varepsilon}(t)^{2}
\Big]\chi_{E_{\varepsilon}}(t)dt+o(\varepsilon^{2}), \quad (\varepsilon\to 0^+).
\end{split}
\end{equation}
\end{proposition}

Now, we can derive the pointwise second-order necessary condition for stochastic optimal controls in the general case. Similar to the convex control constraint case, the solution $y_{1}^{\varepsilon}(\cdot)$ (of the variational equation (\ref{firstvariequ})), which is only an infinitesimal of order $\frac{1}{2}$ with respect to $\varepsilon$ as $\varepsilon\to 0^+$, appears in the variational formulation (\ref{shorttaylor}). Therefore, the Lebesgue differentiation theorem cannot be used directly, either. Hence, the following regularity assumption needs to be introduced.
\begin{enumerate}
\item [{\em (C6)}] For any $v\in V$, $\mss(\cdot,\bar{x}(\cdot),v)\in \ml_{2,\mmf}^{1,2}(\mr)$, and map $v\mapsto \nabla\mss(\t,\bar{x}(\t),v)$ is continuous on $V$ for a.e. $\t\in[0,T]$.
\end{enumerate}

We have the following two results.

\begin{theorem}
\label{2ordercondition nonconvex}
Let (C4)--(C6) hold.  If $\bar{u}(\cdot)$ is a singular optimal control in the sense of Pontryagin-type maximum principle on some control region $V\subset U$, then
\begin{eqnarray*}
& &\mss(\t,\bar{x}(\t), v)
\big(b(\t,\bar{x}(\t),v)-b(\t,\bar{x}(\t),\bar{u}(\t))\big)+ \nabla\mss(\t,\bar{x}(\t),v)
\big(\sigma(\t,\bar{x}(\t),v)-\sigma(\t,\bar{x}(\t),\bar{u}(\t))\big)\nonumber\\
& &
\quad+\frac{1}{2}\mt(\t,\bar{x}(\t),v)
\big(\sigma(\t,\bar{x}(\t),v)-\sigma(\t,\bar{x}(\t),\bar{u}(\t))\big)\nonumber\\
& &\le 0,
\qquad \forall \  v\in V, \ \ a.e.\ (\omega,t)\in \Omega\times[0,T].
\end{eqnarray*}
\end{theorem}

\begin{corollary}\label{2orderconditionth sequalszero}
Let the assumptions of Theorem \ref{2ordercondition nonconvex} hold. If
$$\mss(t,\bar{x}(t), v)=0,\quad \forall v\in V, \quad a.e.\ (\omega,t)\in \Omega\times[0,T],$$
then
\begin{equation*}
\mt(\t,\bar{x}(\t),v)
(\sigma(\t,\bar{x}(\t),v)-\sigma(\t,\bar{x}(\t),\bar{u}(\t)))^2\le 0,
\ \ \ \  \forall \ \ v\in V, \ a.e.\ (\omega,t)\in \Omega\times[0,T].
\end{equation*}
\end{corollary}


\Acknowledgements{This work was partially supported by the National Basic Research Program of China (973 Program) under grant 2011CB808002, the NSF of China under grant 11231007, and the PCSIRT under grant IRT1273 (from the Chinese Education Ministry).}



\begin{thebibliography}{99}
\bahao\baselineskip 11.5pt


\bibitem{Bensoussan81} Bensoussan A.
Lectures on stochastic control, in Nonlinear Filtering and Stochastic Control, 1--62. {Lecture Notes in Math.}, 972, Berlin: Springer-Verlag, 1982

\bibitem{Bismut78} Bismut J -M. An introductory approach to duality in optimal stochastic
control. {SIAM Rev}, 1978, 20: 62--78%

\bibitem{Bonnans12} Bonnans J F, Silva F J. First and second order necessary conditions for stochastic optimal control problems. {Appl. Math. Optim.}, 2012, 65: 403--439%

\bibitem{Gabasov72}  Gabasov R, Kirillova F M. High order necessary conditions for optimality. {SIAM J. Control}, 1972, 10: 127--168%

\bibitem{Haussmann76} Haussmann U G. General necessary conditions for optimal control of stochastic systems. {Math. Prog. Study}, 1976,  6: 34--48%

\bibitem{Kushner72} Kushner H J. Necessary conditions for continuous parameter stochastic
optimization problems. {SIAM J. Control Optim.}, 1972,  10:  550--565%

\bibitem{Nualart06} Nualart D.
The Malliavin Calculus and Related Topics, Second edition. Berlin: Springer-Verlag, 2006%

\bibitem{Peng90} Peng S. A general stochastic maximum principle for optimal control problems. {SIAM J. Control Optim.}, 1990, 28: 966--979%

\bibitem{Tang10} Tang S. A second-order maximum principle for singular optimal stochastic controls. {Discrete Contin. Dyn. Syst. Ser. B}, 2010, 14: 1581--1599%

\bibitem{Yong99} Yong J, Zhou X.
Stochastic Controls: Hamiltonian Systems and HJB Equations. New York: Springer-Verlag, 1999%

\bibitem{zhangH14a} Zhang H, Zhang X. Pointwise second-order necessary conditions for stochastic optimal controls, Part I: The case of convex control constraint. Preprint.%

\bibitem{zhangH14b} Zhang H, Zhang X. Pointwise second-order necessary conditions for stochastic optimal controls, Part II: The general case. Preprint.%

\end{thebibliography}
\end{document}